\magnification=\magstep1   
\input amstex
   \UseAMSsymbols
\input pictex
\hoffset=0truecm \voffset=0truecm 
\NoBlackBoxes
      \font\gross=cmbx10 scaled\magstep2   \font\abs=cmcsc10
   \font\rmk=cmr8   \font\itk=cmti8    \font\ttk=cmtt8

   \def\mod{\operatorname{mod}}
   \def\add{\operatorname{add}}

   \def\End{\operatorname{End}}
   \def\Sing{\operatorname{Sing}}

   \def\ss{\ssize }
\def\arr#1#2{\arrow <1.5mm> [0.25,0.75] from #1 to #2}
 \def\Rahmenbio#1%
   {$$\vbox{\hrule\hbox%
                  {\vrule%
                       \hskip0.5cm%
                            \vbox{\vskip0.3cm\relax%
                               \hbox{$\displaystyle{#1}$}%
                                  \vskip0.3cm}%
                       \hskip0.5cm%
                  \vrule}%
           \hrule}$$}        

\def\Rahmen#1%
   {\centerline{\vbox{\hrule\hbox%
                  {\vrule%
                       \hskip0.5cm%
                            \vbox{\vskip0.3cm\relax%
                               \hbox{{#1}}%
                                  \vskip0.3cm}%
                       \hskip0.5cm%
                  \vrule}%
           \hrule}}}
\headline{\ifnum\pageno=1\hfill %
    \else\ifodd\pageno \hfil\Rechts\hfil \else \hfil\Links\hfil \fi  \fi}
    \def\Links{\abs C\. M\. Ringel}
    \def\Rechts{\abs Gabriel-Roiter inclusions and Auslander-Reiten theory}

\vglue3truecm
\centerline{\gross Gabriel-Roiter inclusions}
	\medskip
\centerline{\gross  and Auslander-Reiten theory.}                    	\bigskip
\centerline{Claus Michael Ringel}     \bigskip

\plainfootnote{}{
{\rmk 2000 \itk Mathematics Subject Classification. \rmk Primary 
16D70, 
16D90, 
16G20, 
16G60, 
16G70, 
18E10. 


}}

	\bigskip
{\baselineskip=9pt
\narrower\narrower{\noindent \rmk ABSTRACT. Let $\ss \Lambda$ be 
an artin algebra. The aim of this paper is to outline a strong relationship
between the Gabriel-Roiter inclusions and the Auslander-Reiten theory.
If $\ss X$ is a Gabriel-Roiter submodule of $\ss Y,$ then 
$\ss Y$ is shown to be a factor module of an indecomposable
module $\ss M$ such that there exists an irreducible monomorphism $\ss X \to M$. 
We also will prove that the monomorphisms in
a homogeneous tube are Gabriel-Roiter inclusions, provided the 
the tube contains a module whose endomorphism ring is a division ring.}\par
}
\bigskip\bigskip

Let $\Lambda$ be an artin algebra, and $\mod \Lambda$ the category of 
$\Lambda$-modules of finite length. 
The basic notion of Auslander-Reiten
theory is that of an irreducible map: these are the maps in the radical of
$\mod \Lambda$ which do not belong to the square of the radical. They are used
in order to define the Auslander-Reiten quiver $\Gamma(\Lambda)$: its
vertices are the isomorphism classes $[X]$ of indecomposable $\Lambda$-modules $X$,
and one draws an arrow $[X] \to [Y]$ provided there exists an irreducible map
$X\to Y$. In 1975,
Auslander and Reiten have shown the existence of Auslander-Reiten sequences:
for any indecomposable non-injective module $X$, there
exists an exact sequence $0 \to X \to Y \to Z \to 0$ such that both maps
$X\to Y$ and $Y \to Z$ are irreducible; here, $Z$ is indecomposable 
(and not projective) and any indecomposable
non-projective module occurs in this way. The sequence is uniquely determined
both by $X$ and by $Z$ and one writes $X = \tau Z,$ and $Z = \tau^{-1}X$ and calls
$\tau$ the Auslander-Reiten translation. The middle term $Y$ is not necessarily
indecomposable: the indecomposable direct summands of $Y'$ are precisely the
modules with an arrow $[X] \to [Y']$, and also precisely the 
modules with an arrow $[Y'] \to [Z]$.  

\bigskip

The {\it Gabriel-Roiter measure} $\mu(M)$ of a $\Lambda$-module $M$
is a rational number  defined inductively as follows: For the zero
module $M = 0$, one sets $\mu(0) = 0.$  If $M\neq 0$ is decomposable,
then $\mu(M)$ is the maximum of $\mu(M')$ where $M'$ is a proper submodule of
$M$, whereas for an indecomposable module $M$, one sets
$$
 \mu(M) = 2^{-|M|}+\max_{M'\subset M} \mu(M').
$$
It is obvious that calculating the maximum of $\mu(M')$, with $M'$ a 
proper submodule of $M$, one may restrict to look at indecomposable
submodules $M'$ of $M$. If $M$ is indecomposable and not simple, then
there always exists an indecomposable submodule $M' \subset M$
such that $\mu(M) - \mu(M') = 2^{-|M|}$, such submodules are called
{\it Gabriel-Roiter submodules} of $M$, and the inclusion map $M' \subset M $
is called a {\it Gabriel-Roiter inclusion.} Note that 
$M$ may have non-isomorphic Gabriel-Roiter submodules, however all
Gabriel-Roiter submodules of $M$ have at least the same length. 
Inductively, we obtain for any indecomposable module $M$ a chain
of indecomposable submodules
$$
 M_1 \subset M_2 \subset \cdots \subset M_{t-1} \subset M_t = M
$$
such that $M_1$ is simple and all the inclusions $M_{i-1} \subset M_i$
for $2 \le i \le t$ are Gabriel-Roiter inclusions, such a sequence is called
a {\it Gabriel-Roiter filtration.} Given such a Gabriel-Roiter filtration,
we have (by definition)
$$
 \mu(M) = \sum\nolimits_{j=1}^t 2^{-|M_j|},
$$
and it will sometimes be convenient to call also the set 
$I = \{|M_1|,\dots,|M_t|\}$  the Gabriel-Roiter
measure of $M$. Thus the Gabriel-Roiter measure $\mu(M)$ of a module $M$ will be 
considered either as a finite set $I$ of 
natural numbers, or else as the rational number $\sum_{i\in I} 2^{-i},$
whatever is more suitable. 
	\bigskip
The paper is divided into two parts. The first part comprises sections 1 to 3; here we
will discuss in which way Gabriel-Roiter inclusions are related 
to Auslander-Reiten sequences. 
	\medskip
{\bf Theorem A.} {\it Assume that $X \subset Y$ is a Gabriel-Roiter inclusion. 
Then there is an irreducible monomorphism $X \to M$ with $M$ indecomposable 
and an epimorphism $M \to Y$ such that the composition $X \to M \to Y$ is injective}
(and therefore also a Gabriel-Roiter inclusion.)
	\medskip
This result may be reformulated as follows: 
{\it Let $X$ be a Gabriel-Roiter submodule of $Y$. Then
there is an irreducible embedding $X \subset M$ and a submodule $U$ of $M$ with
$X\cap U = 0$ such that $M/U$ is isomorphic to $Y$.}
	\medskip
The proof of Theorem A will be given in section 1. Section 2 will exhibit applications,
in particular we will derive some results concerning the existence of indecomposable
submodules of a given module.
Section 3 will use Theorem A in order to
discuss the so-called take-off part of a 
bimodule algebra. 
	\medskip
There will be an intermediate section 4 where we 
will introduce the notion of a piling submodule; this notion will be helpful for the
further discussions. Note that looking at the piling submodules of a module $M$ 
corresponds to the 
process of constructing inductively Gabriel-Roiter filtrations, starting with the simple
submodules of $M$ and going upwards. 
	\bigskip
The second part, sections 5 and 6, deals with 
modules belonging to homogeneous tubes, or, more generally, to modules 
which have a suitable filtration such that all the factors are isomorphic 
to a given indecomposable module $M$.
Recall that a component of the Auslander-Reiten quiver of $\Lambda$ 
is called a {\it homogeneous tube} provided
it is of the form $\Bbb ZA_\infty/\tau.$ The indecomposable modules belonging to 
a homogeneous tube $\Cal T$ will always be labeled as $M[t]$ with $t\in \Bbb N_1$
such that $M[1]$ is of smallest possible length in $\Cal T$ and any $M[t]$
has a filtration with $t$ factors isomorphic to $M[1]$; the module $M[1]$ will
be called the {\it boundary module of} $\Cal T.$ 
Also, a module is said to be a {\it brick} provided its endomorphism ring is a
division ring.
	\medskip
{\bf Theorem B.} {\it Let $\Cal T$ be a homogeneous tube with indecomposable modules
$M[t]$, where $t\in \Bbb N_1$. Let $m \ge 2.$ 
	\smallskip
\item{\rm(a)} Given a Gabriel-Roiter filtration of $M[m]$, then there is 
a submodule of $M[m]$ isomorphic to $M[1]$ which occurs in the filtration.
	\smallskip
\item{\rm(b)} If $M[1]$ is a brick, then $M[m]$ has a unique Gabriel-Roiter submodule,
namely the unique submodule  of $M[m]$ which is isomorphic to $M[m-1]$. \par }
	\medskip
	
The proof of Theorem B will be given in sections 5 and 6. 
 Note that there is a wealth of artin algebras with 
homogeneous tubes: according to Crawley-Boevey, any tame $k$-algebra with $k$ an
algebraically closed field has homogeneous tubes, but there are also many wild artin
algebras having homogeneous tubes. The case of tame hereditary algebras has been
studied very carefully by Bo Chen, in particular see [C, Corollary 4.5] 
which provides a proof of Theorem B in this case.
	\bigskip
The results of this paper have been presented in a series of lectures
in Bahia Bianca (spring 2006), as well as at the university of Bielefeld, 
(2006 and 2007), and the author is grateful for corresponding helpful comments. 
	\bigskip\bigskip
{\bf 1. Proof of Theorem A.}
	\medskip
Here is a more precise statement.
	\medskip
{\bf 1.1. Theorem.} {\it  Assume that $X \subset Y$ is a Gabriel-Roiter inclusion. Let
$f = (f_i)\:X \to \bigoplus_{i=1}^t M_i$ be the source map, with all 
$M_i$ indecomposable. Assume that the maps 
$f_1,\dots, f_r$ are monomorphisms and the remaining ones not. Then there exists
some $i$ with $1  \le i \le r$ and an epimorphism $g\:M_i \to Y$ such that
$gf_i$ is a monomorphism.}
	\medskip
{\bf Second formulation.} {\it Let $X$ be a Gabriel-Roiter submodule of $Y$. Let
$f\:X \subset M$ be a source map for $X$. Decompose $M = \bigoplus_{i=1}^n M_i$ with
$M_i$ indecomposable, and write $f = (f_i)_i$ with $f_i\:X \to M_i$. 
Then there is an index $i$ such that $f_i$ is injective, and there is a submodule
$U\subset M_i$ with $f(X)\cap U = 0$ such that $M_i/U$ is isomorphic to $Y$. }
	\medskip

{\bf Relevance.} Recall that any indecomposable module $Y$ which is not simple
has a Gabriel-Roiter submodule $X$ and $X$ is indecomposable again. Thus, the
theorem asserts that in order to construct all the indecomposable modules $Y$,
one can proceed inductively as follows, starting with  the simple 
modules. In order to find indecomposable modules $Y$ which are not simple,
we consider an indecomposable module $X$ already constructed, 
an irreducible monomorphism 
$f\: X \to M$, and an epimorphism $g\:M \to Y$ such that 
the composition $gf\:X \to Y$ is injective.
Of course, we can assume that $f$ is an embedding. 
The epimorphism $g$
is determined by its kernel $U$, thus by a submodule
$U$ of $M$ such that $X \cap U = 0.$ 
The picture to have in mind is the following:
$$
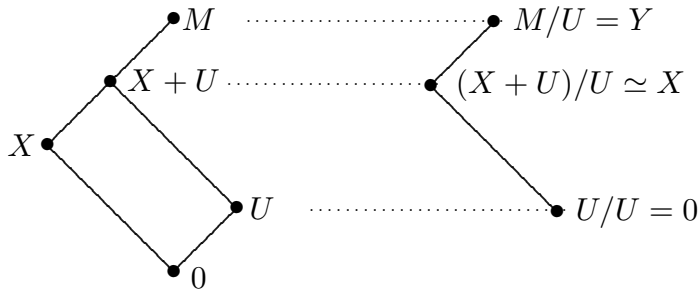
\hbox{\beginpicture
\setcoordinatesystem units <.7cm,.7cm>
\put{\beginpicture
\multiput{} at 0 -.2  4 4.2 /
\plot 2 4  0 2  2 0  3 1  1 3 /
\put{$X$} at -0.4 2
\put{$M$} at 2.4 4
\put{$U$} at 3.4 1
\put{$0$} at 2.4 -0.1
\put{$X+U$} at 2 3
\multiput{$\bullet$} at  2 4  0 2  2 0  3 1  1 3 /
\endpicture} at 0 0 

\put{\beginpicture
\multiput{} at 0 -.2  4 4.2 /
\plot 2 4  1 3  3 1  /

\put{$M/U= Y$} [l] at 2.3 4
\put{$U/U = 0$} [l] at 3.3 1
\put{$(X+U)/U \simeq X$} [l] at 1.4 3 
\multiput{$\bullet$} at  2 4  3 1  1 3 /
\endpicture} at 6 0 

\setdots <1mm>
\plot 1.5 2  5.5 2 /
\plot 1.2 1  4.5 1 /
\plot 2.5 -1  6.5 -1 /
\endpicture}
$$
In this way we obtain all the possible Gabriel-Roiter inclusions, and thus
all the indecomposable modules $Y$. According to 1.2, we only have to look at
finitely many irreducible embeddings $X \to M_i = M$ (and this information is 
stored in the Auslander-Reiten quiver of $\Lambda$). 
The new datum required is the submodule 
$U$ of $M$ with $X\cap U = 0.$ (Unfortunately, we
do not know any criterion on $U$ which tells us whether we obtain
a Gabriel-Roiter inclusion, not even whether we get an indecomposable module $Y$). 

Anyway, we should add that the irreducibility of an embedding $X \to M$ 
yields that for any proper submodule $M'$ of $M$ with $X \subseteq M'$, the
embedding $X \subseteq M'$ splits, thus there is a submodule $U$ of $M$ with
$M' = X\oplus U.$ In this way, the submodules $U$ to be considered correspond to the
proper submodules of $M/X$.

Also, let us stress that here we deal with a quite unusual conjunction of 
indecomposable modules: as we know, the modules $X, M, Y$ are indecomposable.
Since both embeddings of $X$ into $M$ and into $Y$ are mono-irreducible, also
the factor modules $M/X$ and $Y/X \simeq M/(X+U)$  are indecomposable.
	\medskip
{\bf Proof of Theorem 1.1.} 
Let $u\:X \to Y$ be the given Gabriel-Roiter inclusion. 
We denote by $\Sing(X,Y)$ the set 
of maps $f\: X \to Y$ which are not monomorphisms and we know that 
$\Sing(X,Y)$ is closed under addition, see [R4] or [R5].

Let
$f = (f_i)\:X \to \bigoplus_{i=1}^t M_i$ be the source map, with all $M_i$ indecomposable.
We obtain maps $h_i\:M_i \to Y$ such that
$\sum_i h_if_i = u.$ Assume that $h_if_i$ is not a monomorphism,
for $ i > s$, and a monomorphism for the remaining $1\le i\le s$.
If $h_if_i$ is a monomorphism, then also $f_i$ is a monomorphism, thus $s \le r.$

Since $h_if_i$ is not injective for $i > s$, 
also $v = \sum_{ i> s} h_if_i$ is in $\Sing(X,Y)$. Since 
$u$ is a monomorphism, it does not belong to $\Sing(X,Y)$, thus
$u' = u-v=  \sum_{1\le i\le s} h_if_i$ is outside of $\Sing(X,Y)$, 
and therefore a monomorphism (in particular, we see that $s \ge 1).$

Let $Y_i$ be the image of $h_if_i$, for $1 \le i \le s.$ All the following indices are 
$1\le i \le s.$
Assume that none of the maps $h_i$ is surjective. Let $Y_i$ be the image of $h_i$.
Then $\mu(Y_i) < \mu(Y).$
Now under $u'$ the module $X$ embeds into $\bigoplus Y_i$, thus $\mu(X) \le \max \mu(Y_i) <
\mu(Y)$. Since $X$ is a Gabriel-Roiter submodule of $Y$, we must have $\mu(X) = \max\mu(Y_i),$
thus this embedding $X \to \bigoplus M_i \to \bigoplus Y_i$ is a split monomorphism, but then $f$ 
itself is a split monomorphism
--- impossible. This shows that at least one of the $h_i$ is surjective. But then $f_i$ is
a monomorphism (as an irreducible map, it is either mono or epi: if $f_i$ would be epi,
then $|X| > |M_i| \ge |Y|,$ in contrast to the fact that $X$ embeds into $Y$.) 

	\bigskip
This kind of argumentation can be used inductively:
	\medskip
{\bf 1.2. Theorem.} {\it Assume that $u\:X \subset M$ is an 
irreducible monomorphism  with $M$ indecomposable and that  $M \to Y$ 
is an epimorphism such that the composition $X \to M \to Y$ is a Gabriel-Roiter inclusion.
Then 
\item{\rm(a)} either there is a Gabriel-Roiter inclusion of the form
$$
 X @>u>> M = M_0  @>f_1>> M_1  @>>> \cdots  @>f_n>> M_n = Y
$$
with all maps $u, f_1,\dots, f_n$ irreducible (of course $u$ mono), $f_n\cdots f_1$ epi, for some
fixed $n$, or else:

\item{\rm(b)} for every natural number $m$, there is a Gabriel-Roiter inclusion of the form
$$
 X @>u>> M = M_0  @>f_1>> M_1  @>>> \cdots  @>f_m>> M_m @>g>> Y
$$
with all maps $u, f_1,\dots, f_m$ irreducible (of course $u$ mono), $gf_m\cdots f_1$ epi.
\par}
	\medskip
Proof, by induction. The start is given by theorem 1.1: This is just case (b) with
$m = 0.$ Thus, assume that we are in the situation of case (b). Either $g$ is an
isomorphism, then we are in case (a), and nothing else has to be done. If $g$
is not an isomorphism, then it is not split mono (since it is epi), thus we use
the source map for $M_m$ in order to factor $g$. Let $(\phi_i)_i\:M_m \to \bigoplus N_i$
be the source map for $M_m$, we factor $g = \sum \phi_i\psi_i$ with $\psi_i\:N_i\to Y.$

We consider only those indices $i$ such that the composition $\psi_i\phi_i f_1\cdots f_m u$
is mono, say $1 \le  i \le s$.

And we use the usual arguments to see that at least one of the maps  $\psi_i\phi_i f_1\cdots f_m$
with $1 \le i \le s$ has to be surjective, say $i = 1.$ Then let $f_{m+1} = \phi_1$ and use
as new map $g$ the map $\psi_1.$ 
	\bigskip\bigskip
{\bf 2\. Applications.}
	\medskip
An obvious consequence of Theorem A is the following: 
	\medskip
{\bf 2.1. Corollary.} {\it If $X$ is a Gabriel-Roiter
submodule of some module $Y$, then there exists an irreducible monomorphism
$X \to M$ with $M$ indecomposable.}
	\medskip
This shows that a lot of modules cannot be Gabriel-Roiter submodules of other
modules. For example
\item{(1)} Injective modules (of course).
\item{(2)} Let $\Lambda$ be the path algebra of the 
$n$-Kronecker quiver: this is the 
quiver with 2 vertices $a,b$ and $n$ arrows from $a$ to $b$. If $n \ge 2$, then
$\Lambda$ is representation-infinite and has a preinjective component. If $X$ is
an indecomposable preinjective module, then no irreducible map $X \to Y$ with $Y$
indecomposable is a monomorphism.
\item{(3)} Consider the $n$-subspace quiver for some $n\ge 1$
(this is the quiver with $n+1$ vertices, such that there is a unique sink whereas
the remaining vertices are sources, 
and such that there is precisely one arrow from any source to the sink). Let $X$ be an
indecomposable module of the form $\tau^{-t}P$ where $t \ge 1$, where $P$ is the 
unique simple projective module. Then there is
no irreducible map $X \to Y$ with $Y$ indecomposable, which is a monomorphism,
thus $X$ cannot occur as a Gabriel-Roiter submodule. For example, for $n=4$
this concerns all the indecomposable preprojective modules $X$ of length
$6t+1$ with $t\ge 1.$
	\bigskip
Less trivial are the following consequences of Theorem A:
	\medskip
Let $p$ be the maximal length of an indecomposable projective module, 
let $q$ be the maximal length of an indecomposable injective module.
	\medskip

{\bf 2.2. Corollary.} {\it Let $X \to Y$ be a Gabriel-Roiter inclusion. 
Then $|Y| \le pq|X|.$}
	\medskip
Proof: Of course, $X$ cannot be injective. It is well-known that for an 
indecomposable non-injective module $X$, one has $|\tau^{-1}(X)| \le (pq-1)|X|$,
thus the middle term 
$X'$ of the Auslander-Reiten sequence starting in $X$ has length
at most $pq|X|$.  Theorem 1 asserts that $Y$ is a factor module of $X'$, thus also
$|Y| \le pq|X|.$
	\medskip
This result is already mentioned in [R4], as a corollary to Lemma 3.1 of [R4]. Also, 
there we have shown that 2.2 implies the ``successor lemma''. Here are two further
consequences.
	\medskip
{\bf 2.3. Corollary.} {\it Let $M$ be an indecomposable module and $1 \le a < |M|$ a natural
number. Then there exists an indecomposable submodule $M'$ of $M$ with length in the
interval $[a+1,pqa].$}
	\medskip
Proof: Take a Gabriel-Roiter filtration $M_1 \subset \cdots \subset M_n = M.$ Let $i$
be maximal with $|M_i| \le a.$ Then $1\le i <n,$ thus $M_{i+1}$ exists and $a < |M_{i+1}| \le
pq|M_i| \le pqa.$
	\medskip
{\bf 2.4. Corollary.} 
{\it Let $M$ be an indecomposable module and assume that all indecomposable proper
submodules of $M$ are of length at most $b$. Then $|M| \le pqb$.}
	\medskip
Proof: Let $X$ be a Gabriel-Roiter submodule of $M$. By assumption, $|X|\le b,$ thus
$|M| \le pq|X| \le pqb.$
	\medskip
Reformulation: 
Let $\Cal N$ be a class of indecomposable modules. Recall that a module $M$ is said to be 
$\Cal N$-critical provided it does not belong to $\add \Cal N$, but any proper indecomposable
submodule of $M$ belongs to $\Cal N$.
Corollary 2.4 asserts the following: if all the modules in $\Cal N$ are of length at most $b$,
then any $\Cal N$-critical module is of length at most $pqb.$
	\medskip
Observe that the last two corollaries do not refer at all to Gabriel-Roiter notions.
	\bigskip\bigskip

{\bf 3. The take-off part of a bimodule algebra.}
	\medskip
A rational number (or a finite set of natural numbers) will be said to be a
Gabriel-Roiter measure for $\Lambda$, provided there is an indecomposable
$\Lambda$-module with this measure. 
For any Gabriel-Roiter measure $J$, we denote by $\Cal A(J)$ the set of
isomorphism classes of indecomposable modules with measure $J$ (or
representatives of these isomorphism classes).
Recall from [R3] the following: If $\Lambda$ is a representation-infinite
artin algebra, there 
is a countable sequence of Gabriel-Roiter measures $I_1 < I_2 < \cdots$ (the so-called
``take-off measures'') for $\Lambda$ 
such that any other Gabriel-Roiter measure $I$ for $\Lambda$ satisfies $I_t < I$ for
all $t$. Obviously, $I_1 = \{1\}$ and it is easy to see that $I_2 = \{1,t\}$, where
$t$ is the largest possible length of a local $\Lambda$-module of Loewy length $2$. 
	\medskip
In this section, we consider the finite dimensional hereditary algebras with s = 2, where 
$s$ denotes the number of simple modules, thus we deal with representations of
a bimodule ${}_FM_G$, where $F, G$ are division rings.
We assume that $\Lambda = \left[\smallmatrix F & M \cr 0 & G\endsmallmatrix \right]$
is representation-infinite. Of course, we require that $\Lambda$ is an artin algebra, 
thus there is a commutative field $k$ contained in the center both of $F$ and of $G$
and acting centrally on $M$ and such that $\dim_kM$ is finite. 
Let $a = \dim {}_FM, b = \dim M_G.$ The assumption that $\Lambda$ 
is representation-infinite means that $ab \ge 4$.

Often we will present an indecomposable module
by just writing down its dimension vector 
(note that a representation of the bimodule ${}_FM_G$
is a triple $({}_FX,{}_GY,\gamma)$, where $\gamma\: {}_FM_G\otimes_GY \to {}_FX$ is
$F$-linear, its dimension vector is the pair $(\dim {}_FX,\dim {}_GY))$.

Let $P_1, P_2,\dots$ be the sequence of preprojective modules, 
with non-zero maps $P_i \to P_{i+1}.$
$$
\hbox{\beginpicture
\setcoordinatesystem units <1cm,1cm>
\put{$P_1 = (1,0)$} at 0 0
\put{$P_2 = (a,1)$} at 2 1
\put{$P_3 = (ab\!-\!1,b)$} at 4 0
\put{$P_4 = (a^2b\!-\!2a,ab\!-\!1)$} at 6 1
\put{$P_5$} at 8 0

\arr{0.5 0.3}{1.5 0.7}
\arr{2.5 0.7}{3.5 0.3}
\arr{4.5 0.3}{5.5 0.7}
\arr{6.5 0.7}{7.5 0.3}
\setdots<1mm>
\plot 1.2 0  2.6 0 /
\plot 3.2 1  4.2 1 /

\plot 5.4 0  7 0 /
\plot 7.8 1  9 1 /
\plot 8.5 0  9 0 /
\setdots <2mm>
\plot 9 0.5 10 0.5 /

\endpicture}
$$
with $\End(P_{2n-1}) = F,$ and $\End(P_{2n}) = G,$ for all $n$.
	\medskip
Note that always $\Cal A(I_1)$ consists of the simple $\Lambda$-modules.
Here we show that the remaining take-off modules are the modules $P_n$ with $n\ge 2$.
	\medskip
{\bf Proposition.} {\it For $n \ge 2,$ $\Cal A(I_n) = \{P_n\}.$}
	\medskip
For $n=2$, the assertion is true according to the general description of $I_2.$
For $n >2$, we use induction. We have to consider three cases:
	\medskip
{\bf Case 1.} Consider first a bimodule ${}_FM_G$ with  $a, b \ge 2.$
Then all the non-zero maps $P_n \to P_{n+1}$ are monomorphisms. 
Also, since all the irreducible maps ending in $P_n$ are monomorphisms, {\it the monomorphisms
$P_{n-1} \to P_n$ are Gabriel-Roiter inclusions.}

Consider some $n > 2$ and assume that the assertion is true for $n-1$. 
Since there is a Gabriel-Roiter inclusion $P_{n-1} \to P_n$, it follows that $I_{n} =
I_{n-1}\cup\{t\}$ with $t \ge |P_n|.$  
Thus let $Y$ be indecomposable with $\mu(Y) = I_n$, let $X$ be a 
Gabriel-Roiter submodule of $Y$. Then $\mu(Y) = I_{n-1}$, thus by induction $X = P_{n-1}.$
But now we can apply theorem 1.1 above which shows that $Y$ is a factor module of $P_n$.
Since $|Y| = t \ge |P_n|,$ we see that $Y = P_n.$
	\bigskip
{\bf Case 2.} $F \subset G,$ and $M = {}_FG_G,$ thus $a = [G:F].$
Then we deal with the preprojective modules
$$
\hbox{\beginpicture
\setcoordinatesystem units <1cm,1cm>
\put{$P_1 = (1,0)$} at 0 0
\put{$P_2 = (a,1)$} at 2 1
\put{$P_3 = (a\!-\!1,1)$} at 4 0
\put{$P_4 = (a^2\!-\!2a,a\!-\!1)$} at 6 1
\put{$P_5$} at 8 0
\arr{0.5 0.3}{1.5 0.7}
\arr{2.5 0.7}{3.5 0.3}
\arr{4.5 0.3}{5.5 0.7}
\arr{6.5 0.7}{7.5 0.3}
\setdots<1mm>
\plot 1.2 0  2.5 0 /
\plot 3.1 1  4.3 1 /

\plot 5.5 0  6.5 0 /
\plot 7.8 1  9 1 /
\plot 8.5 0  9 0 /

\setdots <2mm>
\plot 9 0.5 10 0.5 /
\endpicture}
$$

\item{$\bullet$} 
{\it The non-zero maps $P_{2n-1} \to P_{2n}$ are injective and are Gabriel-Roiter inclusions.}
\item{$\bullet$} 
The non-zero maps $P_{2n} \to P_{2n+1}$ are surjective.
\item{$\bullet$} 
{\it The non-zero maps $P_{2n-1} \to P_{2n+1}$ are injective and are Gabriel-Roiter inclusions.} 
	\medskip
Consider some $2n$ and assume that the assertion is true for $2n-1$. The argument is the same as
in Case 1, using theorem 1.1.

Also, consider some $2n+1$ and assume that the assertion is true for $2n-1$ and $2n$. Since the 
irreducible maps starting in $P_{2n}$ are epi, we see that $I_{2n+1}$ cannot start with
$I_{2n}.$ Since there are Gabriel-Roiter inclusions $P_{2n-1} \to P_{2n+1},$ we see that
$I_{2n+1} = I_{2n-1}\cup\{t\}$ with $|P_{2n}| > t \ge |P{2n+1}|$. 

Thus let $Y$ be indecomposable with $\mu(Y) = I_{2n+1}$, let $X$ be a 
Gabriel-Roiter submodule of $Y$. Then $\mu(Y) = I_{2n-1}$, thus by induction $X = P_{2n-1}.$
But now we can apply 1.2. 
It shows that $Y$ is a factor module of $P_{2n1}$.
Since $|Y| = t \ge |P_{2n+1}|,$ we see that $Y = P_{2n+1}.$
	\medskip
{\bf Case 3.} $G \subset F,$ and $M = {}_FF_G,$ thus $b = [F:G].$
Then we deal with the preprojectives
$$
\hbox{\beginpicture
\setcoordinatesystem units <1cm,1cm>
\put{$P_1 = (1,0)$} at 0 0
\put{$P_2 = (1,1)$} at 2 1
\put{$P_3 = (b\!-\!1,b)$} at 4 0
\put{$P_4 = (b\!-\!2,b\!-\!1)$} at 6 1
\put{$P_5$} at 8 0

\arr{0.5 0.3}{1.5 0.7}
\arr{2.5 0.7}{3.5 0.3}
\arr{4.5 0.3}{5.5 0.7}
\arr{6.5 0.7}{7.5 0.3}
\setdots<1mm>
\plot 1.2 0  2.6 0 /
\plot 3.2 1  4.5 1 /

\plot 5.4 0  7 0 /
\plot 7.8 1  9 1 /
\plot 8.5 0  9 0 /

\setdots <2mm>
\plot 9 0.5 10 0.5 /

\endpicture}
$$

The non-zero maps $P_{2n-1} \to P_{2n}$ are surjective, for $n\ge 2$, whereas $P_1 \to P_2$
is injective (and this is a Gabriel-Roiter inclusion).

\item{$\bullet$} 
{\it The non-zero maps $P_{2n} \to P_{2n+1}$ are injective and are Gabriel-Roiter inclusions.}
\item{$\bullet$} 
{\it The non-zero maps $P_{2n} \to P_{2n+2}$ are injective and are Gabriel-Roiter inclusions.} 
	\medskip
Proof: As in case 2, but taking into account
the additional Gabriel-Roiter inclusion $P_1 \to P_2$.
	
	\bigskip\bigskip
{\bf 4. Piling submodules.}
	\medskip
We call an indecomposable submodule $U$ of some module $Y$ {\it piling}, provided
$\mu(V) \le \mu(U)$ for all indecomposable submodules $V$ of $Y$ with $|V| \le |U|$;
actually, it is sufficient to check the condition for the indecomposable submodules 
$V$ of $Y$ with $|V| < |U|$. Namely, if there exists an indecomposable submodule
$V$ of $Y$ with $|V| = |U|$ and $\mu(V) > \mu(U)$, then there exists a proper
submodule $V'$ of $V$ with $\mu(V') > \mu(U)$). 
There is the following alternative description: 
	\medskip
{\bf 4.1. Lemma.} {\it An indecomposable submodule $U$
of $V$ is piling if and only if $\mu(Y)$ starts with $\mu(U)$}
(this means that  $\mu(U) = \mu(Y)\;\cap\;\{1,2,\dots,|U|\}$).
	\medskip
Proof: Let $U$ be an indecomposable submodule of $Y$, let $U_1 \subset U_2 \subset 
\cdots \subset U_s = U$ and $Y_1 \subset Y_2 \subset \cdots \subset Y_t$ be 
Gabriel-Roiter filtrations, where $Y_t \subseteq Y$ is an indecomposable submodule
with $\mu(Y_t) = \mu(Y).$ 

First, assume that $U$ is piling in $Y$. We claim that $|U_i| = |Y_i|$ for 
$1 \le i \le s$. If not, then there is some minimal $i$ with $|U_i| \neq |Y_i|$
and since $\mu(U) \le \mu(Y)$, we must have $|U_i| > |Y_i|.$ But then $V = Y_i$
is a submodule of $Y$ with $|V| < |U|$ and $\mu(V) = \{Y_1,\dots,Y_i\} >
\mu(U),$ a contradiction to the assumption that $U$ is piling in $Y$.

Second, assume that $\mu(U) = \mu(Y)\;\cap\;\{1,2,\dots,|U|\}.$
Let $V$ be an indecomposable submodule of $Y$ with $|V| \le |U|$
and assume that $\mu(V) > \mu(U)$.
If $V_1 \subset V_2 \subset \cdots \subset V_r = V$ is a Gabriel-Roiter
filtration of $V$, then there must be some $1 \le j \le \min(s+1,r)$ such that
$|U_i| = |V_i|$ for $1 \le i < j$ and either $j = s+1$ 
or else $|U_j| > |V_j|.$ The case $j = s+1$ cannot happen, since otherwise
$|U| = |U_s| = |V_s| < |V_{s+1}| \le |V|,$ but $|V| \le |U|.$
Thus we have $|V_j|< |U_j| .$ But then $\mu(V) > \mu(Y)$, since
$|V_i| = |U_i| = |Y_i|$ for $1 \le i < j$ and $|V_j| < |U_j| = |Y_j|$.
This is impossible: a submodule $V$ of $Y$ always satisfies $\mu(V) \le \mu(Y).$
	\medskip
Note that all the submodules in a Gabriel-Roiter
filtration of an indecomposable module are piling, but usually there are additional
ones: for example all simple submodules are piling. 
The fact that a submodule $U$ of $Y$ is piling depends only on the isomorphism class
of $U$ and the set of isomorphism classes of submodules $V$ of $Y$
with $|V| \le |U|$ (but for example not on the embedding of $U$ into $Y$).
Here are some further properties:
	\medskip
{\bf 4.2.}
{\it Assume that $U \subseteq V \subseteq W$. If $U$ is piling in $V$ and $V$ is piling in $W$, then $U$ is piling in $W$.} Proof: 
If $\mu(W)$ starts with $\mu(V)$ and $\mu(V)$ starts with $\mu(U)$, then
obviously $\mu(W)$ starts with $\mu(U).$
	\medskip
{\bf 4.3.} {\it If $U \subseteq V \subseteq W$ and $U$ is a piling submodule of $W$, then also of $V$.} Proof: Let $X$ be a submodule of $V$ with $|X| \le |U|$. Consider $X$ as a
submodule of $W$ and conclude that $\mu(X) \le \mu(U).$
	\medskip
{\bf 4.4.} {\it If $U \subseteq X\oplus Y$ is a piling submodule, then at least one of the
maps $U \to X$ or $U \to Y$ is an embedding with piling image.} 
Recall the strong Gabriel property: Assume that $U, X_1,\dots, X_n$ are
indecomposable modules and there are given maps $f_i\:U \to X_i$ such that the
map $f = (f_i)_i\:U \to \bigoplus X_i = X$ is a monomorphism and its image is a
piling submodule of $X$. Then at least one of the maps $f_i$ is a monomorphism
(and its image is a piling submodule of $X_i$).
Now let $U \subseteq X\oplus Y$ be a piling submodule.
According to the strong Gabriel property, one of the maps $U \to X$, $U \to Y$
is an embedding, say $f\:U \to X$. Since $U$ is piling in $X\oplus Y,$ it follows
that $f(U)$ is piling in $X$.
	\bigskip\bigskip
{\bf 5. Piling submodules of modules with a homogeneous $M$-filtration.}
	\medskip
{\bf 5.1.}
Let $M$ be an indecomposable module. An {\it $M$-filtration}
of a module $Y$ is a chain of modules
$$
 M = M[1] \subset M[2] \subset \cdots \subset M[m] = Y
$$
such that $M[i]/M[i-1] \simeq M$ for all $2 \le i \le m;$ 
the number $m$ is called the {\it $M$-length} of $Y$. Observe that a module $Y$ with an  $M$-filtration is just an iterated self-extension of the module $M$.

In case all the
modules $M[i]$ are indecomposable and the inclusion maps are mono-irreducible maps,
we call this a {\it homogeneous} $M$-filtration.
For example, if $\Cal T$ is a homogeneous tube 
of the Auslander-Reiten quiver of $\Lambda$ and $M$ is the boundary module of $\Cal T,$
then any module belonging to $\Cal T$ has a homogeneous $M$-filtration.

On the other hand, it is easy to construct modules $M$ with self-extensions
$$
 0 \to M \to M[2] \to M \to 0
$$
where $M[2]$ is indecomposable
such that the inclusion $M \to M[2]$ is not mono-irreducible. For example, take the quiver
with vertices $a,b,c$, one arrow $a \to b$, two arrows $b \to c$. Then
any indecomposable module $M$ with dimension vector $(1,1,1)$ is as required.
	\medskip
The aim of this section is to study Gabriel-Roiter filtrations of 
modules with a homogeneous $M$-filtration. 
	\bigskip
{\bf 5.2. Theorem.} {\it Let $Y$ be a module with a homogeneous $M$-filtration.
Let $U$ be a piling submodule of $Y$.}

(a) {\it If $|U| \le |M|,$ then $U$ is isomorphic to a submodule of $M$.}

(b) {\it If $|U| \ge |M|$, 
then any Gabriel-Roiter filtration of $U$ contains a module isomorphic to $M$.}
	\medskip
{\bf 5.3. Corollary.} {\it Let $Y$ be a module with a homogeneous $M$-filtration.
Then any Gabriel-Roiter filtration of $Y$ contains a module isomorphic to $M$.}
	\medskip
This is the special case of (b) where $U = Y.$ On the other hand, the assertion (a)
of Theorem B is a direct application of 5.3: Given a homogeneous tube $\Cal T$ with
modules $M[t]$, then any module $M[t]$ has a homogeneous $M[1]$-filtration.
	\bigskip
Before we start with the proof of Theorem 5.2, we insert a general observation:
	\medskip
{\bf 5.4.} {\it Let $X \to Y$ be a Gabriel-Roiter inclusion.
Assume that $Y'$ is a submodule of $Y$ with 
$X+Y' = Y.$ Then $|Y'| > |Y/X|.$}
	\medskip
{\bf Proof.} The submodule 
$Y'$ of $Y$ maps onto
$Y'/(X\cap Y') \simeq (X+Y')/X = Y/X$ with kernel $X\cap Y'.$
If $|Y'| \le |Y/X|,$ then this surjective
map has to be an isomorphism, thus $X\cap Y' = 0.$ But then $Y = X\oplus Y',$
whereas $Y$ is indecomposable and $0 \neq X \neq Y.$
	\bigskip
Proof of 5.2(a). Assume that $|U| \le |M|$. According to 5.4, we see that
$M[m-1]+U$ is a proper submodule of $Y$, thus $M[m-1]+U = M[m-1]\oplus U'$
for some submodule $U'$ of $Y$ which is isomorphic to a proper submodule of $M$.
Now $U$ is a piling submodule of $Y$, thus also of $M[m-1]\oplus U'$,
therefore of $M[m-1]$ or of $U'$. In the first case, use induction on $m$. 
In the second case, just recall that $U'$ is isomorphic to a submodule of $M$.
	\medskip
Proof of 5.2(b). We can assume that $M$ is not simple, since otherwise $M[m]$
is serial and nothing has to be shown. 

Let $U_1 \subset U_2 \subset \cdots \subset U_s = U$
be a Gabriel-Roiter filtration of $U$. Assume that 
$|U_r| < |M|$ and $|U_{r+1}| \ge |M|$
by assumption, such an $r$ must exist, since $M$ is not simple and $|U_s| \ge |M|.$
We apply (a) to the submodule $U_r$ (as a piling submodule of $U$ it is piling in $Y$)
and see that $U_r$ is isomorphic to a submodule $M'$ of $M$. From the definition of 
a Gabriel-Roiter filtration it follows that $U_{r+1} \le |M|.$ Thus
$|U_{r+1}| = |M|.$ Now we apply (a) to the sold submodule $U_{r+1}$ and see
that $U_{r+1}$ is isomorphic to $M$. 
	\bigskip\bigskip
{\bf 6. Modules with homogeneous $M$-filtrations, where $M$ is a brick.}
	\medskip
{\bf 6.1.}
We will assume now in addition that the endomorphism ring of $M$ is a division ring,
thus that $M$ is a brick.  Then we can use the process of simplification [R1]:
Let $\Cal F(M)$ be the full subcategory of all modules which have an $M$-filtration.
The new assumption implies that this category $\Cal F(M)$ 
is an abelian category, even a length category, and $M$ is its only simple object.
Of course, the $M$-filtrations of an object $Y$ are just the composition series
of $Y$ when considered as an object of $\Cal F(M).$
Thus, if $Y$ has an $M$-filtration
$$
 M = M[1] \subset M[2] \subset \cdots \subset M[m] = Y
$$
with all $M[i]$ for $1\le i \le m$ indecomposable, then $Y$ considered as an
object of $\Cal F(M)$ is uniform, thus $M$ is the only submodule of $Y$ isomorphic
to $M$. 
	\bigskip
In case $Y$ has a unique $M$-filtration
$$
 M = M[1] \subset M[2] \subset \cdots \subset M[m] = Y
$$
then $Y$, considered as an object in $\Cal F(M)$ is even serial and then
all the factors $M[t]/M[s]$ with $0 \le s < t \le m$ are indecomposable
(here, $M[0] = 0$). 

Conversely, if $Y$ has the $M$-filtration
$$
 M = M[1] \subset M[2] \subset \cdots \subset M[m] = Y
$$
and all the factor modules $M[i]/M[i-2]$ with $2 \le i \le m$ are indecomposable
(again, we set $M[0] = 0$), then $Y$ has only one $M$-filtration. 
	\bigskip
{\bf 6.2. Theorem.} {\it Assume that $M$ is a brick.
Let $Y$ be a module with a homogeneous $M$-filtration
and assume that
$$
 M = M[1] \subset M[2] \subset \cdots \subset M[m] = Y
$$
is the only $M$-filtration of $Y$. If $U$ is a piling submodule of $Y$ and 
$|U| \ge |M|$, then $U = M[j]$ for some $1 \le j\le m$.}
	\medskip
Remark: The assumption $|U| \ge |M|$ is important. For example, if $|M| \ge 2,$
then the socle of $M[m]$ is of length at least $m$, thus for $m \ge 2$, 
there are many simple submodules and all are piling. 
	\medskip
Proof of 6.2. We can assume that $M$ is not simple, since otherwise $M[m]$ is 
a serial module and the submodules $M[j]$ are the only non-zero submodules.

We use induction on $m$. For $m = 1$ nothing has to be shown. Thus let $m \ge 2.$.
Let $U$ be a piling submodule of $Y$ and $|U| \ge |M|$. We use a second induction,
now on $|U|$ in order to show that $U = M[j]$ for some $j$. 
The induction starts with $|U| = |M|.$ In this case theorem 1.1 shows that
$U$ is isomorphic to $M$, but according to the process of simplification, 
any submodule of $Y$ isomorphic to $M$ is equal to $M$, thus $U = M = M[1]$.

Now let $|U| > |M|.$ Let $U'$ be a Gabriel-Roiter
submodule of $U$. 

First, assume that $|U'| < |M|.$ Then, according to 5.2(a), there is a submodule $M'$
of $M$ which is isomorphic to $U'$. Using the definition of a Gabriel-Roiter
sequence, and the fact that there is the inclusion $M' \subset M$ with $M$
indecomposable, we see that $|U| \le |M|.$ Together with $|U| \ge |M|$, it follows
that $|U| = |M|.$ Using again 5.2(a), now for $U$, we see that $U$ is isomorphic to
$M$. However, according to the process of simplification, 
any submodule of $Y$ isomorphic to $M$ is equal to $M$, thus $U = M = M[1]$.

Next, consider the case $|U'| \ge |M|.$ By the second induction, we see that
$U' = M[j]$ for some $j\ge 1$. Of course, since $U'$ is a proper submodule of $Y$,
we see that $j < m$. The inclusion $M[j] \subset M[j+1]$ and the definition
of a Gabriel-Roiter filtration shows that $|U| \le |M[j+1]|.$
Claim: We can assume that $M[m-1]+U$ is a proper submodule of $M[m]$. 

Otherwise, we easily see
that $M[m]/M[j] = M[m-1]/M[j] \oplus U/M[j]$. But the process of simplification
shows that $M[m]/N[j]$ is indecomposable. 
Also, $U/M[j] = U/U'$ is non-zero,
thus $M[m-1]/M[j] = 0$ and $U = M[m].$ 

Since $M[m-1]+U$ is a proper submodule of $M[m]$, we see that
$M[m-1]+U = M[m-1]\oplus C$ for some submodule $C$ of $Y$ and $C$ is
isomorphic to a proper submodule of $M$. 
Write the inclusion map $U \to M[m-1]\oplus C$ in the form $[f,f']^t,$ where
$f\:U \to M[m-1]$ and $f' \to C.$ Since $U$ is a piling submodule of $Y$
and $U \subseteq M[m-1]\oplus C \subset Y$, it follows that $U$
is a piling submodule of $M[m-1]\oplus C$, thus either $f$ or $f'$ is an
embedding with piling image. However $|U| > |M|$ whereas $|C| < |M|$, thus
$f\:U \to M[m-1]$ is an embedding, and its image is a piling submodule. 
Now we use the induction on $m$ in order to conclude that $f(U) = M[i]$ for
some $i$. In particular, $U$ is isomorphic to $M[i]$. 
But there is no non-zero homomorphism $M[i] \to C$, since $C$ is a proper
submodule of $M$ (using simplification). This shows that $f' = 0$ and therefore
the embedding $U \to M[m-1]\oplus C$ is just the map $f\:U \to M[m-1]$:
this shows that $U=M[i].$ 
	\bigskip
{\bf 6.3. Corollary.} {\it Assume that $M$ is a brick.
Let $Y$ be a module with a homogeneous $M$-filtration
and assume that
$$
 M = M[1] \subset M[2] \subset \cdots \subset M[m] = Y
$$
is the only $M$-filtration of $Y$. If $m \ge 2$, 
then $Y$ has precisely one Gabriel-Roiter submodule, 
namely $M[m-1].$}
	\medskip
Proof of 6.3. Let $X$ be a Gabriel-Roiter submodule of $Y = M[m]$,
where $m \ge 2.$  
A Gabriel-Roiter submodule is a piling submodule. If $|X| \le |M|$, then by 5.2(a),
$X$ is isomorphic to a submodule $M'$ of $M$. With $X$ also $M'$ is a Gabriel-Roiter
submodule of $Y$, thus the inclusions $M' \subseteq M \subset M[2]
\subseteq M[m]$ show that
$M' = M$ (and $m = 2).$ 
Thus $X$ is isomorphic to $M$. However the process
of simplification asserts that $M$ is the only submodule of $M[m]$ isomorphic to
$M$. This shows that $X = M[1]$.

Thus we can assume that $|X| > |M|$ and use theorem 6.2.
	\bigskip
{\bf 6.4.} Proof of Theorem B (b). 
Let $\Cal T$ 
be a homogeneous tube of the Auslander-Reiten quiver of $\Lambda$.
Denote the modules in $\Cal T$ by $M[m]$, where
$|M[m]| = m|M|$ and $M = M[1]$. Assume in addition that $M$ is a brick. 
Since $M$ is a brick,
we can consider the abelian category $\Cal F(M)$ which contains all the
modules $M[m]$. We have a chain of inclusions
$$
  M = M[1] \subset M[2] \subset \cdots
$$
and for $m \ge 2$, the factor module $M[m]/M[m-2]$ is indecomposable
(here again, we set $M[0] = 0).$
This shows that any $M[m]$ has the unique $M$-filtration
$$
  M = M[1] \subset M[2] \subset \cdots \subset M[m],
$$
thus we can use the previous corollary.
	\bigskip
{\bf 6.5. Remarks.} In 6.3 as well as in Theorem B (b), 
the assumption $m \ge 2$ is important: the module $M = M[1]$ usually
will have more than one Gabriel-Roiter submodules. For example consider the
four-subspace-quiver and $\Cal T$ a homogeneous tube containing a module $M$ of length 6.
The module $M$ has 4 maximal submodules and all are Gabriel-Roiter submodules.
	\medskip
Well-known examples of homogeneous tubes such that 
the endomorphism ring of the boundary module $M$ is a division ring
are the homogeneous tubes of a tame hereditary 
algebra, of a tubular algebra or of a canonical algebra [R2]. As we have mentioned
in the introduction, tame hereditary algebras have been considered by Bo Chen in [C].
	\medskip
For a tubular algebra, the boundary modules of homogeneous tubes are of
unbounded length. The same is true in case $\Lambda$ is a tame hereditary or
a canonical $k$-algebra and 
the algebraic closure of $k$ is not a finite field extension
of $k$.
	\bigskip
At the end of the paper, let us consider also modules of infinite lengths.
We consider again a homogeneous tube $\Cal T$ with indecomposable
modules $M[m]$, where $t\in \Bbb N_1.$ There is a chain of irreducible maps
$$
 M[1] \to M[2] \to \cdots \to M[m] \to M[m+1] \to \cdots,
$$
and we denote by $M[\infty]$ the corresponding direct limit; such a module is
called a {\it Pr\"ufer module.} 
	\medskip
{\bf 6.6. Corollary.} {\it Let $\Cal T$ be a homogeneous tube with indecomposable
modules $M[m]$ and assume that $M[1]$ is a brick. 
Then the Gabriel-Roiter measure of $M[\infty]$ is a rational number.}
	\medskip
Proof. The module $M[\infty]$ has a Gabriel-Roiter filtration starting with a
Gabriel-Roiter filtration of $M[1]$ and then using precisely the modules $M[m]$.
Let $X$ be a Gabriel-Roiter submodule of $M[1]$. Note that the length of
$M[m]$ is $sm$, with $s = |M[1]|.$ 
Thus
$$
 \gamma(M[\infty]) = \gamma(X) + \sum_{m\ge 1} 2^{-|M[m]|} =
     \gamma(X) + \sum_{m\ge 1} 2^{-sm} = \gamma(X) + \frac1{2^s-1}.
$$
	\bigskip\bigskip
{\bf References.}\medskip
{\frenchspacing
\item{[C]} Chen, B.:
  Comparison of Auslander-Reiten theory and Gabriel-Roiter measure 
  approach to the module categories of tame hereditary algebras. 
  Comm. Algebra 36(2008), 4186-4200 
\item{[R1]} Ringel, C. M.: Representations of K-species and bimodules. J. Algebra 41 (1976), 269-302. 
\item{[R2]} Ringel, C. M.: {\it 
Tame algebras and integral quadratic forms.} Springer LNM 1099 (1984). 
\item{[R3]} Ringel, C. M.: The Gabriel-Roiter measure.
 Bull. Sci. math. 129 (2005), 726-748.
\item{[R4]} Ringel, C. M.: Foundation of the Representation Theory of Artin Algebras, 
 Using the Gabriel-Roiter Measure. 
 Contemporary Math. 406. Amer.Math.Soc. (2006), 105-135. 
\item{[R5]} Ringel, C. M.: The theorem of Bo Chen and Hall polynomials. 
 Nagoya Journal 183 (2006).
\par}
	\bigskip\bigskip

{\rmk Fakult\"at f\"ur Mathematik, Universit\"at Bielefeld,
POBox 100\,131, \ 
D-33\,501 Bielefeld}

{\rmk E-mail address:} {\ttk ringel\@math.uni-bielefeld.de}

\bye